\documentclass[12pt]{amsart}
\usepackage[margin=1in]{geometry}
\usepackage{times}
\usepackage{amsthm,amssymb}
\usepackage{amsmath,xypic}
\usepackage[usenames,dvipsnames]{color}
\usepackage[backend=biber, style=alphabetic, citestyle=authoryear-comp, firstinits=true, maxcitenames=50]{biblatex}
\usepackage[
colorlinks=true,hyperindex, citecolor=cyan, urlcolor=cyan]{hyperref}
\usepackage{epigraph}
\usepackage{graphicx}
\usepackage{fancyhdr} 
\usepackage{enumitem}
\newcommand{\be}{\begin{equation}}
\newcommand{\ee}{\end{equation}}
\newcommand{\bes}{\begin{equation*}}
\newcommand{\ees}{\end{equation*}}
\newcommand{\bea}{\begin{eqnarray}}
\newcommand{\eea}{\end{eqnarray}}
\newcommand{\beas}{\begin{eqnarray}}
\newcommand{\eeas}{\end{eqnarray}}
\newcommand{\ben}{\begin{note}}
\newcommand{\een}{\end{note}}
\newcommand{\bexl}{\vskip0.1em\noindent\hrulefill\vskip1em\begin{ExerciseList}}
\newcommand{\eexl}{\end{ExerciseList}\hrulefill}

\newcommand{\bthm}{\begin{theorem}}
\newcommand{\ethm}{\end{theorem}}
\newcommand{\bpro}{\begin{prop}}
\newcommand{\epro}{\end{prop}}
\newcommand{\bcor}{\begin{corollary}}
\newcommand{\ecor}{\end{corollary}}
\newcommand{\bcon}{\begin{conjecture}}
\newcommand{\econ}{\end{conjecture}}
\newcommand{\bp}{\begin{proof}}
\newcommand{\ep}{\end{proof}}
\newcommand{\blem}{\begin{lemma}}
\newcommand{\elem}{\end{lemma}}
\newcommand{\bn}{\begin{note}}
\newcommand{\en}{\end{note}}
\newcommand{\benum}{\begin{enumerate}}
\newcommand{\eenum}{\end{enumerate}}
\newcommand{\bed}{\begin{defn}}
\newcommand{\eed}{\end{defn}}
\newcommand{\brem}{\begin{remark}}
\newcommand{\erem}{\end{remark}}

\newcommand{\btik}{\begin{tikzpicture}\begin{axis}[scale=0.5,axis y line=center, axis x line=middle]}
\newcommand{\etik}{\end{axis}\end{tikzpicture}}

\let\into=\hookrightarrow

\newcommand{\upperRomannumeral}[1]{\uppercase\expandafter{\romannumeral#1}}

\newtheorem{theorem}[equation]{Theorem}      
\newtheorem{lemma}[equation]{Lemma}          %
\newtheorem{corollary}[equation]{Corollary}  
\newtheorem{proposition}[equation]{Proposition}

\theoremstyle{definition}

\theoremstyle{definition}
\newtheorem{defn}[equation]{Definition}
\theoremstyle{remark}

\theoremstyle{definition}
\newtheorem{remark}[equation]{Remark}

\numberwithin{equation}{section}

\let\into=\hookrightarrow
\let\isom=\simeq

\let\tensor=\otimes

\newcommand{\abs}[1]{\left\vert#1\right\vert}

\newcommand{\Pic}{{\rm Pic}}

\newcommand{\Z}{{\mathbb Z}}

\renewcommand{\int}{\operatorname{int}}
\renewcommand{\O}{{\mathcal O}}
\renewcommand{\P}{{\mathbb P}}

\newcommand{\mapright}[1]{{\xymatrix{{}\ar[r]^{#1}&{}}}}

\renewcommand{\bpro}{\begin{proposition}}
\renewcommand{\epro}{\end{proposition}}

\let\citep=\cite 

\begin{document}

\title[]{Corrigendum: On the construction of Weakly Ulrich bundles}%
\author{Kirti Joshi}%
\address{Math. department, University of Arizona, 617 N Santa Rita, Tucson
85721-0089, USA.} \email{kirti@math.arizona.edu}

\thanks{}%

\let\cite=\parencite

\begin{abstract}
This note corrects a mistake in   \cite[Theorem~4.1]{joshi2021a}.  \emph{The error noted here  does not affect any other results of \cite{joshi2021a}.}  To correct the error, here I  prove a more general result (Theorem~\ref{th:corrigendum}) and deduce Theorem~\ref{th:corrigendum2-gen} and also the correct version of the previously announced theorem  in Theorem~\ref{th:corrigendum2}. Theorem~\ref{th:pic-Z-pluricanonical-case} supplements Theorem~\ref{th:corrigendum2-gen}. In it, I prove that if  $k$ is an algebraically closed field of characteristic $p\geq 3$ and $X/k$ is any smooth, projective, minimal surface  of general type and with $\Pic(X)=\Z$, then for all integers $r\geq 5$,  $X$ is embedded as a smooth surface by its pluricanonical linear system $X\into \abs{\omega_X^r}$, and $E=F_*(\omega_X^{r+1})(1)$ is an almost Ulrich bundle  for the pluricanonical embedding $X\into\abs{\omega_X^r}$ of $X$ and for the ample line bundle provided by this embedding. Corollary~\ref{cor:corrigendum} generalizes \cite[Theorem 3.1]{joshi2021a}. 
\end{abstract}
\maketitle
%

\newcommand{\pn}[1]{\mathbb{P}^{#1}}
\newcommand{\pnn}{\pn{n}}
\newcommand{\pnno}{\O_{\pnn}}
\section{Introduction} This note corrects an error in \cite[Theorem 4.1]{joshi2021a}. 
Unfortunately because of this error, \cite[Theorem 4.1]{joshi2021a} is not correct in that the claimed bundle $E=F_*(\omega_X)$ fails to be an almost Ulrich bundle, but  in the correction presented below (Theorem~\ref{th:corrigendum2}) I show (with no new hypothesis) that $E=F_*(\omega^2_X)$ is an almost Ulrich bundle. \emph{This error does not affect any of the other theorems proved in \cite{joshi2021a}.} I first prove a more general result (see Theorem~\ref{th:corrigendum}) and deduce Theorem~\ref{th:corrigendum2-gen} and also the correct version  (Theorem~\ref{th:corrigendum2}) of the previously announced result \cite[Theorem 4.1]{joshi2021a}    as a consequence.  Theorem~\ref{th:pic-Z-pluricanonical-case} supplements Theorem~\ref{th:corrigendum2}--it shows that  provided $p\geq 3$ and $r\geq 5$, any smooth, projective  surface of general type with $\Pic(X)=\Z$ is embedded as a smooth surface by its pluricanonical embedding $X\into\abs{\omega_X^r}$ and $F_*(\omega_X^r)(1)$ is an almost Ulrich bundle on $X$. I also note that Corollary~\ref{cor:corrigendum} provides examples of Theorem~\ref{th:corrigendum} and also generalizes \cite[Theorem 3.1]{joshi2021a}. In particular surfaces satisfying hypotheses of Theorem~\ref{th:corrigendum}, Theorem~\ref{th:corrigendum2-gen} and Theorem~\ref{th:corrigendum2} exist (Remarks~\ref{rem:1} and \ref{rem:2}). 

\section{Error located in \cite[Theorem~4.1]{joshi2021a}}
Before providing the correct version of the Theorem and its proof, let me point out the error  in \cite[Theorem~4.1]{joshi2021a}. I will use the notation of \cite{joshi2021a} for the discussion of the error.

The error occurs in the following step of the proof of \cite[Theorem~4.1]{joshi2021a}: it is claimed in the proof that $H^2(E(-m))=0$ for all $m\leq 1$. As $E=F_*(\omega_X)(1)$, this says
$$H^2(E(-m))=H^2(F_*(\omega_X)\tensor\O_X(1-m))=0.$$ \emph{Unfortunately this is false for $m=1$}. Indeed for $m=1$ this is
$$H^2(E(-1))=H^2(F_*(\omega_X))=H^2(\omega_X)\isom H^0(\O_X)\neq 0.$$ [I had claimed in my proof that the vanishing of $H^2(F_*(\omega_X))$ follows from that of $H^2(\omega^p_X)=H^0(\omega_X^{1-p})=0$. This  vanishing is true by ampleness of  $\omega_X$ assumed in \cite[Theorem~4.1]{joshi2021a}, but the equality $H^2(\omega_X)=H^2(\omega_X^p)$ claimed in my proof is false as is clear from  this discussion.]

\section{Corrected Theorem}
Now let me prove a general result (Theorem~\ref{th:corrigendum}) and deduce from it  Theorem~\ref{th:corrigendum2} which replaces \cite[Theorem~4.1]{joshi2021a}. Let me also mention that the result given below includes as a special case \cite[Theorem~3.1]{joshi2021a} and also provides new examples of almost Ulrich bundles not covered by \cite[Theorem~3.1]{joshi2021a}.  For   the definition of Ulrich, weakly Ulrich and almost Ulrich bundles see \cite[Section 2]{joshi2021a}.

\bthm\label{th:corrigendum}
Let $k$ be an algebraically closed field of characteristic $p>0$, $X/k$ be a smooth, projective surface equipped with a projective embedding $(X\into \P^n,\O_X(1))$. Let $$F:X\to X$$ be the absolute Frobenius morphism of $X$. Let $M$ be a line bundle on $X$ and $\omega_X$ the canonical line bundle of $X$.
Assume that the following hold
\benum[label={{\bf(\arabic{*})}}]
\item $M$ is an ample line bundle on $X$;
\item  $\omega_X\tensor M$ is an arithmetically Cohen-Macaulay (ACM) line bundle on $X$ i.e. 
$$H^1(\omega_X\tensor M\tensor\O_X(m))=0\quad \text{for all }m\in\Z.$$
\item $H^0(\omega_X\tensor M\tensor\O_X(-p))=0$.
\eenum
Then $E=F_*(\omega_X\tensor M)(1)$ is an almost Ulrich bundle on $X$.
\ethm

\brem 
It is important to note that Theorem~\ref{th:corrigendum} does not assume that Kodaira vanishing holds for $X$.
\erem

\bp
 To prove that $E$ is an \emph{almost Ulrich bundle}, I have to prove vanishing of cohomology $H^i(E(-m))=0$ for $i=0,1,2$ and for $m\in\Z$ in the following ranges:
 \benum[label={{\bf(V.\arabic{*})}},start=0]
   \item\label{v0} $H^0(E(-m))=0$ for all  $m\geq 2$,
 \item\label{v1} $H^1(E(m))=0$ for all $m\in\Z$ (i.e. $E$ is an ACM bundle)
  \item\label{v2} $H^2(E(-m))=0$ for all $m\leq 1$.
 \eenum

The following  standard facts are used in the proof and proofs of these facts are readily found in \cite{hartshorne-algebraic}. 
\benum[label={{\bf(\arabic{*})}}]
\item as $X$ is smooth and projective,  $F$ is a finite flat morphism and hence for any vector bundle $V,M$ on $X$, one has 
$$H^i(F_*(V))=H^i(V).$$
\item The  projection formula for $F$: $$F_*(V)\tensor M=F_*(V\tensor F^*(M)).$$ 
\item Finally, if $M$ is a line bundle on $X$ then $$F^*(M)=M^{\tensor p}.$$
\eenum

Let me now prove the vanishing assertions \ref{v0}, \ref{v1} and \ref{v2}.

First let me  prove \ref{v1}. I claim that $\omega_X\tensor M$ is ACM implies \ref{v1} i.e. $H^1(E(m))=0$ for all $m\in \Z$ i.e. $E$ is also an ACM bundle. Indeed using the projection formula and the facts listed above one  has  
$$F_*(\omega_X\tensor M)(1-m)=F_*\left(\omega_X\tensor M\tensor F^*(\O_X(1-m))\right)=F_*\left(\omega_X\tensor M\tensor \O_X(p(1-m))\right)$$ 
and hence
$$H^1(E(m))=H^1(F_*(\omega_X\tensor M)(1-m))=H^1(\omega_X\tensor M\tensor \O_X(p(1-m)))=0$$
and as  I have assumed $H^1(\omega_X\tensor M\tensor \O_X(m))=0$ for all $m\in\Z$, the required vanishing follows. 
So this settles the vanishing required \ref{v1} for $H^1$.

The next step is  to prove \ref{v2} i.e. $H^2(E(-m))=0$ for $m\leq 1$. Observe again that
$$H^2(E(-m))=H^2(F_*(\omega_X\tensor M)(1-m))=H^2(\omega_X\tensor M\tensor\O_X(p(1-m)))$$
Note that for any $m\leq 1$ one has $1-m\geq 0$ and as $M$ is an ample line bundle, so for all $m\leq 1$, the line bundle $M\tensor\O_X(p(1-m))$ is ample and so its dual $M^{-1}\tensor\O_X(-p(1-m))$ is an anti-ample line bundle. By Serre duality one has
$$ H^2(\omega_X\tensor M\tensor\O_X(p(1-m)))= H^0(M^{-1}\tensor\O_X(-p(1-m)))=0 \text{ by anti-ampleness.}$$
So this proves the required vanishing \ref{v2} for $H^2$.

So it remains to prove \ref{v0} i.e. $H^0(E(-m))=0$ for $m\geq 2$. Again
$$H^0(E(-m))=H^0(F_*(\omega_X\tensor M)(1-m))=H^0(\omega_X\tensor M\tensor \O_X(p(1-m))).$$
As $m\geq 2$, one has $1-m\leq -1$ and so $p(1-m)\leq -p$ for all $m\geq 2$. This means that one has 
inclusion of ideal sheaves (for $m\geq 2$) $\O_X(p(1-m)) \subseteq \O_X(-p)$ and hence
a short exact sequence
$$0\to \O_X(p(1-m)) \to \O_X(-p)\to \O_X(-p)/\O_X(p(1-m))\to 0.$$
Tensoring this by $\omega_X\tensor M$ and noting that $\omega_X\tensor M$ is locally free, one has an exact sequence (for all $m\geq 2$)
$$0\to \omega_X\tensor M\tensor\O_X(p(1-m)) \to \omega_X\tensor M\tensor\O_X(-p)\to \omega_X\tensor M\tensor\left(\O_X(-p)/\O_X(p(1-m))\to 0\right).$$
Taking cohomology  one has
$$0\to H^0(\omega_X\tensor M\tensor \O_X(p(1-m))) \to H^0( \omega_X\tensor M\tensor \O_X(-p)) \to \cdots$$
and so the required vanishing follows from the vanishing of $H^0(\omega_X\tensor M\tensor \O_X(-p))=0$ which is my hypothesis.
\ep

Now let me record a non-trivial consequence of Theorem~\ref{th:corrigendum} from which the correction to \cite[Theorem 4.1]{joshi2021a} will be deduced as a consequence in Theorem~\ref{th:corrigendum2}.

\bthm\label{th:corrigendum2-gen}
Suppose $X/k$ is a smooth, projective,  surface over an algebraically closed field of characteristic $p>0$ and equipped with a projective embedding $X\hookrightarrow \P^n$ with $\O_X(1)$ provided by this embedding. Assume
\benum[label={{\bf(\arabic{*})}}] 
\item $X$ is minimal,
\item $\omega_X^r=\O_X(1)$  for some integer $r\geq 1$ (so $\omega_X$ is ample, and hence $X$ is minimal and of general type)
\item $H^1(X,\O_X)=0$,
\item $p\geq 3$.
\eenum
Then $E=F_*(\omega_X^{r+1})(1)$ is an almost Ulrich bundle on $X$.
\ethm
\bp 
This is a consequence of Theorem~\ref{th:corrigendum}.   To invoke Theorem~\ref{th:corrigendum}, one needs to verify all the hypotheses of that assertion are satisfied. 

First take $$M=\omega_X^r=\O_X(1),$$ so $M$ is ample. This choice of $M$ gives $$E=F_*(\omega_X\tensor M)(1)=F_*(\omega_X^{r+1})(1).$$  I claim that with this choice of $M$ all the hypotheses of Theorem~\ref{th:corrigendum} are satisfied and hence $E=F_*(\omega_X^{r+1})(1)$  is an almost Ulrich bundle.

Let me begin by proving that $\omega_X\tensor M$ is an ACM line bundle. This means I have to show that $H^1(\omega_X\tensor M\tensor \O_X(m))=0$ for all $m\in\Z$. But
$$H^1(\omega_X\tensor M\tensor \O_X(m)) = H^1(\omega_X^{r+1}\tensor\O_X(m))$$
holds for any $m\in\Z$. 
Since $\omega_X^r=\O_X(1)$, so one has  that $$\omega_X^{r+1}\tensor\O_X(m)=\omega_X^{r+1+rm}.$$

Hence to prove the vanishing of $H^1(\omega_X\tensor M\tensor \O_X(m))=0$ for all $m\in\Z$ it is enough to prove that
$$H^1(\omega_X^\ell)=0$$
for all $\ell\in\Z$.

If $\ell=0$ then by definition  $\omega_X^0=\O_X$ and so one wants to prove $H^1(\O_X)=0$, but this is one of the hypothesis of Theorem~\ref{th:corrigendum2}. If $\ell=1$ then $H^1(\omega_X)=H^1(\O_X)$ by Serre duality and again this is zero by hypothesis. If $\ell\geq 2$ then $H^1(\omega_X^\ell)=H^1(\omega_X^{1-\ell})$ and as $\ell\geq 2$ one has $1-\ell \leq -1$. 

So the claim that $H^1(\omega_X^\ell)=0$ for all $\ell\in\Z$ is reduced to the claim that $H^1(\omega_X^{\ell})=0$ for all $\ell\leq -1$. If Kodaira Vanishing Theorem holds for $X$ this is immediate from ampleness of $\omega_X$ for ${\ell}\leq -1$. However I do not assume Kodaira Vanishing Theorem holds for $X$.

The way around unavailability of Kodaira Vanishing Theorem is to use an equally delicate result of \cite[II, Theorem 1.7]{ekedahl88}. This requires that $X$ is minimal of general type and $p\geq3$. This is the case here as  $p\geq 3$,  $X$ is minimal and also  of general type as $\omega_X=\O_X(1)$ is ample by my assumptions. Thus one deduces from \cite[II, Theorem 1.7]{ekedahl88} that for $p\geq 3$ one has $H^1(\omega_X^\ell)=0$ for all $\ell\leq -1$ and hence I have proved that $H^1(\omega_X^\ell)=0$ for all $\ell\in\Z$. Thus $\omega_X\tensor M$ is an ACM bundle which proves that the second hypothesis of Theorem~\ref{th:corrigendum} holds in the present situation.

Finally to check the remaining assumption of Theorem~\ref{th:corrigendum} is also available under the hypothesis of Theorem~\ref{th:corrigendum2} note that
$H^0(\omega_X\tensor M\tensor \O_X(-p))=H^0(\omega^{r+1}_X\tensor\omega_X^{-rp})$ as $\omega_X^r=\O_X(1)=M$. If $p\geq3$ then this cohomology is certainly zero as $\omega_X^{r+1-rp}$ is anti-ample for any $p\geq 3$ and any $r\geq 1$ as $r+1-rp\leq r+1-3r=1-2r\leq -1$ for $r\geq 1$ (note that I have again used the hypothesis that $p\geq 3$).
\ep

Now let me record, in Theorem~\ref{th:corrigendum2}, the correct version of \cite[Theorem 4.1]{joshi2021a}. This corrected version   is a special case of Theorem~\ref{th:corrigendum2-gen} obtained  by taking $r=1$. Note that the hypothesis of Theorem~\ref{th:corrigendum2} below are the same as that of \cite[Theorem 4.1]{joshi2021a} but itemized here for additional clarity.

\bthm\label{th:corrigendum2}
Suppose $X/k$ is a smooth, projective,  surface over an algebraically closed field of characteristic $p>0$ and equipped with a projective embedding $X\hookrightarrow \P^n$ with $\O_X(1)$ provided by this embedding. Assume
\benum[label={{\bf(\arabic{*})}}] 
\item $X$ is minimal,
\item $\omega_X=\O_X(1)$   (so $\omega_X$ is ample, and hence $X$ is minimal and of general type)
\item $H^1(X,\O_X)=0$,
\item $p\geq 3$.
\eenum
Then $E=F_*(\omega_X^{2})(1)$ is an almost Ulrich bundle on $X$.
\ethm

\bp 
This is immediate from Theorem~\ref{th:corrigendum2-gen} upon taking $r=1$.
\ep

The following theorem provides a large class of examples of surfaces satisfying the hypothesis of Theorem~\ref{th:corrigendum2-gen}.

\bthm\label{th:pic-Z-pluricanonical-case}
Let $X/k$ be a smooth, projective surface over an algebraically closed field of characteristic $p>0$. Assume the following hold:
\benum[label={{\bf(\arabic{*})}}] 
\item $p\geq 3$,
\item $\Pic(X)\isom \Z$,
\item $X$ is of Kodaira dimension two.
\eenum
Then for every integer $r\geq 5$ one has
\benum[label={{\bf(\arabic{*})}}]
\item  $X$ is embedded as a  smooth surface $$\varphi_r:X\into \abs{\omega_X^r}=\P(H^0(X,\omega_X^r))$$ by its pluricanonical linear system $\abs{\omega_X^r}$, and
\item  $E=F_*(\omega_X^{r+1})(1)$ is an almost Ulrich bundle on $X$ for its $r$-canonical projective embedding $\varphi_r:X\into \abs{\omega_X^r}$.
\eenum
\ethm

\bp 
Let me note that my hypothesis $\Pic(X)\isom \Z$ says that the Picard scheme of $X$ is reduced  and its connected component is equal to zero and hence its tangent space $H^1(X,\O_X)=0$. Let $H$ be a generator for $\Pic(X)$. One can assume that $H$ is the class of an ample line bundle as $X$ is projective. As $X$ is of general type, $\omega_X$ is a big divisor and hence $\omega_X$ is a positive multiple of $H$ and so $\omega_X$ is ample. One also sees from this that $X$ is minimal, indeed, as every irreducible curve  $C\subset X$ is a multiple of the ample class $H$ and hence in particular for every irreducible curve $C$ on $X$, one has $C\cdot C\geq 1$. So $X$ contains no exceptional curves and hence $X$ is minimal.  Now consider the pluricanonical linear system  $\abs{\omega_X^r}$ and the pluricanonical (rational) mapping
$$\varphi_r:X \to \abs{\omega_X^r}.$$
As $X$ is of Kodaira dimension two
one has $\dim(\varphi_r(X))=2$ for all sufficiently large $r$. As $X$ is minimal, $p\geq 3$ and $X$ is general type, by \cite[Main Theorem]{ekedahl88}, it is enough to take $r\geq 5$ instead of $r$ sufficiently large for $\varphi_r$ to be a morphism. By \cite{mumford62} (which is an appendix to \cite{zariski62}), the morphism $$\varphi_r:X\to \abs{\omega_X^r}$$ factors as 
$$ \varphi_r:X\to X^{can}\into \abs{\omega_X^r},$$
where $X^{can}={\rm Proj}(\oplus_{m\in\Z}H^0(\omega_X^m))$ is the canonical model of $X$ and the morphism $X\to X^{can}$ is birational, and contracts all the (finitely many)  smooth rational curves with self-intersection $-2$ in $X$ and is an isomorphism outside its exceptional locus. 

As I have remarked earlier, under my hypothesis $\Pic(X)=\Z$,  every irreducible curve $C$ lying on $X$ is a multiple  of the ample class $H$. Hence one has $C\cdot C\geq 1$ and so there are no  curves  on $X$ with $C\cdot C=-2$. Thus the exceptional locus of $X\to X^{can}$ is empty and hence $X\to X^{can}$ is an isomorphism and $X\isom X^{can}\into \abs{\omega_X^r}$. So $X$ is pluricanonically embedded by $\varphi_r$ as a smooth surface  for every $r\geq 5$. 

Now one can apply Theorem~\ref{th:corrigendum2-gen} to $X\into \abs{\omega_X^r}$ and deduce that  $E=F_*(\omega^{r+1}_X)$ is an almost Ulrich bundle on $\varphi_r:X\mapright{\isom} X^{can}\subset \abs{\omega_X^r}$. This completes the proof.
\ep

The following corollary of Theorem~\ref{th:corrigendum} illustrates that this theorem is more general than \cite[Theorem~3.1]{joshi2021a} as it provides new examples of bundles even for surfaces in $\P^3$.

\begin{corollary}\label{cor:corrigendum}
Let $X\subset\P^3$ be a smooth surface of degree $d$. Assume
\benum[label={{\bf(\arabic{*})}}]
\item $M=\O_X(r)$ with $r\geq 1$, and
\item  $d-4+r<p$.
\eenum 
 Then $E=F_*(\omega_X(r))(1)$ is an almost Ulrich bundle on $X$. In particular taking $M=\O_X(1)$, one obtains the case considered in \cite[Theorem 3.1]{joshi2021a}.
\end{corollary}
\bp 
Since $X\subset \P^3$ is a smooth surface of degree $d\geq 5$, one has $\omega_X=\O_X(d-4)$. As $M=\O_X(r)$ with $r\geq 1$ so $M$ is ample.  Also one has $\omega_X\tensor M=\O_X(d-4+r)$ and so $\omega_X\tensor M=\O_X(d-4+r)$ is clearly ACM line bundle on $X$ (as $H^1(\O_X(\ell))=0$ for any $\ell\in\Z$ and for any smooth surface in $\P^3$). Moreover $\omega_X\tensor M\tensor\O_X(-p)=\O_X(d-4+r-p)$. So the last hypothesis of Theorem~\ref{th:corrigendum} i.e. $H^0(\omega_X\tensor M\tensor\O_X(-p))=0$ is satisfied if one assumes condition $d-4+r-p<0$. So the corollary follows.

Now taking $M=\O_X(1)$ (i.e. $r=1$), and hence $\omega_X\tensor M=\O_X(d-3)$ and then $\omega_X\tensor M$ is clearly ACM line bundle on $X$. The last hypothesis is the condition $d-4+1-p=d-3-p<0$. But this is precisely the situation of \cite[Theorem~3.1]{joshi2021a}. Thus \cite[Theorem~3.1]{joshi2021a} is a special case of the above Theorem~\ref{th:corrigendum}.
\ep

\begin{remark}\label{rem:1}
In particular  Corollary~\ref{cor:corrigendum} also shows that surfaces $X$ and line bundles $M$ satisfying all of the hypothesis of Theorem~\ref{th:corrigendum} exist and so the Theorem~\ref{th:corrigendum} is non-vacuous. 
\end{remark} 	

\begin{remark}\label{rem:2}
	In the notation of Corollary~\ref{cor:corrigendum}, let  $X\subset\P^3$ be a smooth quintic hypersurface i.e. $d=5$ and take $r=1$ in Corollary~\ref{cor:corrigendum}, so $$\omega_X=\O_X(d-4)=\O_X(1)=M$$
	and by Corollary~\ref{cor:corrigendum} one sees that $E=F_*(\omega_X^2)(1)$ is a weakly Ulrich bundle on $X$. This provides examples of surfaces satisfying Theorem~\ref{th:corrigendum2-gen} and Theorem~\ref{th:corrigendum2}. 
\end{remark}

\printbibliography
\end{document}